\documentclass[12pt]{article}
\usepackage{amsfonts}

\newcommand{\Aut}{\mathop{\mathrm{Aut}}}
\newcommand{\BA}{\mathop{\mathrm{BAut}}}
\newcommand{\bU}{\mathbb{U}}
\newcommand{\bQ}{\mathbb{Q}}
\newcommand{\bZ}{\mathbb{Z}}
\newtheorem{theorem}{Theorem}
\newtheorem{prop}[theorem]{Proposition}
\newtheorem{cor}[theorem]{Corollary}
\newtheorem{lemma}[theorem]{Lemma}
\newtheorem{unnumber}{}

\newtheorem{prepf}[unnumber]{Proof}
\newenvironment{pf}{\prepf\rm}{\endprepf}
\newcommand{\qed}{\qquad$\Box$}

\begin{document}
\title{Some isometry groups of Urysohn space}
\author{P. J. Cameron and A. M. Vershik\\[10mm]
\small School of Mathematical Sciences, Queen Mary, University of London,\\
\small Mile End Road, London E1 4NS, U.K.\\[5mm]
\small St.~Petersburg Department of Steklov Institute of Mathematics,\\
\small Russian Academy of Sciences,\\
\small Fontanka 27, St.~Petersburg 191023, Russia}
\date{}
\maketitle

\begin{abstract}
We construct various isometry groups of Urysohn space (the unique complete
separable metric space that is universal and homogeneous), including
abelian groups which act transitively, and free groups which are dense in
the full isometry group.
\end{abstract}

\section{Introduction}

In a posthumously-published paper, P.~S.~Urysohn~\cite{urysohn} constructed
a remarkable complete separable metric space $\bU$ that is both
\emph{homogeneous} (any isometry between finite subsets of $\bU$ can be
extended to an isometry of $\bU$) and \emph{universal} (every complete
separable metric space can be embedded in $\bU$). This space is unique up
to isometry.

The second author~\cite{vershik} showed that $\bU$ is both the generic complete
metric space with distinguished countable dense subset (in the sense of
Baire category) and the random such space (with respect to any of a wide
class of measures).

In this paper, we investigate the isometry group $\Aut(\bU)$ of $\bU$, and
construct a few interesting subgroups of this group.

Our main tool is an analogous countable metric space $\bQ\bU$, the
unique countable homogeneous metric space with rational distances.
The existence and uniqueness of $\bQ\bU$ follow from the
arguments used to establish the existence and uniqueness of $\bU$
in~\cite{urysohn}. Alternatively, this can be deduced from the
fact that the class of finite metric spaces with rational
distances has the amalgamation property, which, by
Fra\"{\i}ss\'e's theorem~\cite{fraisse}, implies the result. Now
$\bU$ is the completion of $\bQ\bU$ (see more details in
\cite{vershik}). In particular, any isometry of $\bQ\bU$ extends
uniquely to~$\bU$. Our notation suggests that $\bQ\bU$ is ``rational
Urysohn space''.

Let $\Aut(\bQ\bU)$ and $\Aut(\bU)$ be the isometry groups of $\bQ\bU$
and $\bU$. We show that $\Aut(\bQ\bU)$ as a subgroup of $\Aut(\bU)$ is
dense in $\Aut(\bU)$ in the natural topology induced by the
product topology on $\bU^\bU$). We also show that $\bQ\bU$ has an
isometry which permutes all its points in a single cycle (indeed,
it has $2^{\aleph_0}$ conjugacy classes of such isometries). The
closure of the cyclic group generated by such an isometry is an
abelian group which acts transitively on $\bU$, {\it so that
$\Aut(\bU)$ is  a monothetic group and $\bU$ carries an abelian
group structure (in fact, many such structures)}. Moreover, the
free group of countable rank acts as a group of isometries of
$\bQ\bU$ which is dense in the full isometry group (and hence is also
dense in $\Aut(\bU)$).

The universal rational metric space $\bQ\bU$ is characterized
by the following property: If $A,B$ are finite metric spaces with
rational distances (we say \emph{rational metric spaces} for
short) with $A\subseteq B$, then any embedding of $A$ in $\bQ\bU$
can be extended to an embedding of $B$. It is enough to assume
this in the case where $|B|=|A|+1$, in which case it takes a more
convenient form:
\begin{description}
\item{($*$)} If $A$ is a finite subset of $\bQ\bU$ and $g$ is a function from $A$ to
the rationals satisfying
  \begin{itemize}
  \item $g(a)\ge 0$ for all $a\in A$,
  \item $|g(a)-g(b)|\le d(a,b)\le g(a)+g(b)$ for all $a,b\in A$,
  \end{itemize}
then there is a point $z\in \bQ\bU$ such that $d(z,a)=g(a)$ for all $a\in A$.
\end{description}

Furthermore, $\bQ\bU$ is homogeneous (any isometry between finite
subsets of~$\bQ\bU$ extends to an isometry of $\bQ\bU$), and every
countable rational metric space can be embedded isometrically in
$\bQ\bU$.

Note that the same definition and above condition of
universality of metric spaces is valid if instead of the field
$\Bbb Q$ or $\Bbb R$ we consider any countable additive subgroup
of $\Bbb R$, for example, the group of integers; in this case we have the integer
universal metric space ${\Bbb Z}\bU$, which  was considered in
\cite{hco}; we will use it in Section 4; another example ---
the universal metric space with possible values of metric
$\{0,1,2\}$ --- is simply the universal graph $\Gamma$
(see~\cite{Ca,vershik1}) with the ``edge'' metric on the set of its
vertices.

\section{$\Aut(\bQ\bU)$ is dense in $\Aut(\bU)$}

The weak topology on the group $\Aut(\bU)$ of isometries of $\bU$ is that induced by the product
topology on $\bU^\bU$. In particular, $g_n\to g$ if and only if, for any finite sequence
$(u_1,\ldots,u_m)$ of points and any $\epsilon>0$, there exists $n_0$ such that
$d(g_n(u_i),g(u_i))<\epsilon$ for $1\le i\le m$ and $n\ge n_0$.
\begin{theorem}
The group $\Aut(\bQ\bU)$ is a dense subgroup of $\Aut(\bU)$ in weak topology.
\end{theorem}
It suffices to show the following property of $\bQ\bU$:

\begin{prop}
Given $\epsilon>0$ and $v_1,\ldots,v_n,v_1',\ldots,v_{n-1}',v_n''\in \bQ\bU$
such that $(v_1,\ldots, v_{n-1})$ and $(v_1',\ldots,v_{n-1}')$ are isometric
and
\[|d(v_i',v_n'')-d(v_i,v_n)|<\epsilon,\]
there exists $v_n'\in \bQ\bU$ such that $(v_1,\ldots,v_n)$ and
$(v_1',\ldots,v_n')$ are isometric and  $d(v_n',v_n'')<\epsilon$.
\end{prop}
\begin{pf}
{Assuming this for a moment, we complete the proof of the density as follows. We are given an
isometry $g$ of $\bU$ and points $u_1,\ldots,u_m\in \bU$. Choose points $v_1,\ldots,v_m\in \bQ\bU$
with $d(v_i,u_i)<\epsilon/4m$. Now using  the above proposition, we inductively choose points
$v_1',\ldots,v_m'$ so that $(v_1,\ldots,v_m)$ and $(v_1',\ldots,v_m')$ are isometric and
$d(v_i',g(u_i))<i\epsilon/m$. For suppose that $v_1',\ldots, v_{n-1}'$ have been  chosen. Choose
any point $v_n''\in \bQ\bU$ with $d(g(u_n),v_n'')<\epsilon/4m$. Then
\[d(u_i,u_n)-\epsilon/2m<d(v_i,v_n)<d(u_i,u_n)+\epsilon/2m,\]
and
\[d(g(u_i),g(u_n))-(4i+1)\epsilon/4m<d(v_i',v_n'')
<d(g(u_i),g(u_n))+(4i+1)\epsilon/4m,\]
so
\[|d(v_i',v_n'')-d(v_i,v_n)|<(4i+3)\epsilon/4m\le(4n-1)\epsilon/4m.\]
So we may apply the proposition to choose $v_n'$ with
$d(v_n',v_n'')<(4n-1)\epsilon/4m$. Then $d(v_n',g(u_n))<n\epsilon/m$, and
we have finished the inductive step. At the conclusion, we have
$d(v_n',g(u_n))<n\epsilon/m\le\epsilon$ for $1\le n\le m$.

Now we find an isometry of $\bQ\bU$ mapping $v_i$ to $v_i'$ for $1\le i\le m$ (by the homogeneity
of~$\bQ\bU$), and the proof is complete.}\qed
\end{pf}

\paragraph{Proof of Proposition 1.} We have to extend the set
$\{v_1',\ldots,v_{n-1}',v_n''\}$ by adding a point $v_n'$ with prescribed
distances to $v_1',\ldots,v_{n-1}'$ and distance less than $\epsilon$ to
$v_n''$. So it is enough to show that these requirements do not conflict,
that is, that
\[|d(v_n',x)-d(v_n',y)|\le d(x,y)\le d(v_n\,x)+d(v_n\,y)\]
for $x,y\in\{v_1',\ldots,v_{n-1}',v_n''\}$. There are no conflicts if
$x,y\ne v_n''$: this follows from the fact that the points $v_1,\ldots,v_n$
exist having the required distances. So we may assume that $x=v_i'$ and
$y=v_n''$, in which case the consistency follows from the hypothesis.\qed

\section{$\BA(\bU)$ is dense in $\Aut(\bU)$}

For a metric space $M$, we define $\BA(M)$ to be the group of all
\emph{bounded} isometries of $M$ (those satisfying $d(x,g(x))\le k$ for all
$x\in M$, where $k$ is a constant). Clearly it is a normal subgroup of
$\Aut(M)$, though in general it may be trivial, or it may be the whole
of $\Aut(M)$.

We show that $\BA(\bQ\bU)$ is a dense subgroup of $\Aut(\bQ\bU)$: in other words,
any isometry between finite subsets of $\bQ\bU$ can be extended to a bounded
isometry of $\bQ\bU$. This is immediate from the following lemma.

\begin{lemma}
Let $f$ be an isometry between finite subsets $A$ and $B$ of $\bQ\bU$, satisfying
$d(a,f(a))\le k$ for all $a\in A$. Then $f$ can be extended to an isometry
$g$ of $\bQ\bU$ satisfying $d(x,g(x))\le k$ for all $x\in \bQ\bU$.
\label{bddense}
\end{lemma}

\begin{pf}
Suppose that $f:a_i\mapsto b_i$ for $i=1,\ldots,n$, with $d(a_i,b_i)\le k$.
It is enough to show that, for any point $u\in \bQ\bU$, there exists $v\in \bQ\bU$
such that $d(b_i,v)=d(a_i,u)$ for all $i$ and $d(u,v)\le k$. For then we can
extend $f$ to any further point; the same result in reverse shows that we can
extend $f^{-1}$, and then we can construct $g$ by a back-and-forth argument.

The pont $v$ must satisfy $d(b_i,v)=d(a_i,u)$ and $d(u,v)\le k$. We must show
that these requirements are consistent; then the existence of $v$ follows
from the extension property of $\bQ\bU$. Clearly the consistency conditions for
the values $d(b_i,v)$ are satisfied. So the only possible conflict can arise
from the inequality
\[|d(v,u)-d(v,b_i)|\le d(u,b_i)\le d(v,u)+d(v,b_i).\]
We wish to impose an upper bound on $d(v,u)$, so a conflict could arise only
if a lower bound arising from the displayed equation were greater than $k$,
that is, $|d(v,b_i)-d(u,b_i)|>k$, or equivalently, $|d(u,a_i)-d(u,b_i)|>k$.
But this is not the case, since
\[|d(u,a_i)-d(u,b_i)|\le d(a_i,b_i)\le k.\]\qed\end{pf}

\section{A cyclic isometry of universal spaces $\bU$.}

Universal Urysohn metric space has very important and rather
surprising property which we formulate in the following theorem:

\begin{theorem}
{There is an isometry $g \in ISO(\bU)$ of the universal Urysohn
space $\bU$ such that the $\langle g\rangle$-orbit  of some point
$x$ (i.e., the set $\{g^n x; n\in \Bbb Z$\}) is dense in $\bU$; in
this case the $\langle g\rangle$-orbit of each point is dense in
$\bU$.}
\end{theorem}

The second claim follows from the first one directly. One of the
important corollaries of this theorem is

\begin{prop}
 There are transitive abelian groups of isometries of $\bU$ of
 infinite exponent.
\end{prop}

\begin{pf}{Let $\overline{G}$ be the closure in $\Aut(\bU)$
 of the cyclic group $\langle g\rangle$. Since the orbits of $g$
 are dense, it is clear that $\overline{G}$ is transitive.
 Moreover, as the closure of an abelian group, it is itself
 abelian. For, if $h,k\in\overline{G}$, say $h_i\to h$ and $k_i\to
 k$; then $h_ik_i=k_ih_i\to hk=kh$.}\qed
\end{pf}

 What is the structure of the group $\overline{G}$ --- the closure of
 the group ${\Bbb Z}=\{g^n,n\in \Bbb Z\}$?
 Since there are many choices for such $g$, we must expect that their
 closures will not all be alike. In particular, there should be
 some choices of $g$ such that $\overline{G}=G(g)$ is torsion-free,
 and others for which it is not. What is very important is the fact
 that the Urysohn space can be equipped with the structure of an abelian
 group; this must help to find some appropriate  model of this
 space; we will discuss this question elsewhere.

 The proof of Theorem 4 follows easily from the analogous fact for
 the rational case:

\begin{theorem}
 {There is an isometry $g$ of $\bQ\bU$ such that $\langle g\rangle$
 is transitive on~$\bQ\bU$.}
\end{theorem}

 Indeed, since the completion of the space $\bQ\bU$ is the universal
 metric space $\bU$, a transitive isometry of $\bQ\bU$ extends to an
 isometry  of $\bU$ with dense orbit.

 We can also put the same question about the universal homogeneous
 {\it integer} metric space $\Bbb Z \bU$;  the answer is as follows:

\begin{theorem}[see \cite{hco}]
 There exists a transitive isometry of the universal integer
 metric space $\Bbb Z\bU$.
\end{theorem}

As we will see, Theorem 6 on the rational case is a corollary of Theorem 7 on the existence of
universal integer metric space. The latter fact was discovered in~\cite{hco}; to make this paper
self-contained, and for its own interest, below we give the proof, which is similar to the
considerations from~\cite{hco}.

We must define a metric with integer, rational, or real values on
the set of integers $\Bbb Z$ that is shift-invariant. We will
call it a cyclic metric. Such a metric is completely determined
by the function $f(i)=d(i,0)$ on the non-negative integers; for
$d(i,j)=f(|j-i|)$. The function should satisfy the constraints

\begin{description}
\item{(a)} $f(i)\ge0$, with equality if and only if $i=0$.
\item{(b)} $|f(i)-f(j)|\le f(i+j)\le f(i)+f(j)$ for all $i,j$.
\end{description}

 We call a function satisfying (a) and (b) a \emph{Toeplitz
 distance function}, and denote the set of such functions
 by $\bQ T$ (resp. $\bZ T$) if the values of $f$ are rational
 (resp. integer).
 If $i,j$ in (a), (b) run over $1, \dots, n$ only, then
 the set of such functions will be denoted by $\bQ T_n$ ($\bZ T_n$).

 Now the cyclic metric space given by such a function is isometric
 to the universal space $\bQ\bU$ if and only if $f$ has the following property:
 \begin{description}
 \item{(c)} given any function $h$ from $\{1,\ldots,k\}$ to the positive
 rationals satisfying
 \[|h(i)-h(j)|\le f(|i-j|)\le h(i)+h(j)\]
 for $i,j\in\{1,\ldots,k\}$, there exists a positive integer $N$ such that
 $h(i)=f(N-i)$ for all $i\in\{1,\ldots,k\}$.
 \end{description}

 We say that a Toeplitz distance function is \emph{universal} if it satisfies (c).
 The same criterion of universality is true for metrics with integer values.
 It is convenient in the proof to use integer metrics  instead of
 rational ones.

 Given $f\in \bQ T_n$ (or $\bZ T_n$), we say that an $m$-tuple $(h(1),\ldots,h(m))$ is
\emph{$f$-admissible} if
 \[|h(i)-h(i+k)|\le f(k)\le h(i)+h(i+k)\]
 for $1\le i< i+k\le m$ and $k\le n$. This notion of admissibility
 agrees with the admissibility of vectors with respect to the
 distance matrix $\{d(i,j)\}_{i,j=1}^n$ in the sense of
 \cite{vershik}. Note that if $h$ is $f$-admissible, then it is
 admissible with respect to the restriction of $f$ to
 $\{1,\ldots,n'\}$ for any $n'\le n$.

 In order to prove the existence of a cyclic isometry of the universal
 rational or integer metric space, it is enough to prove the following theorem.

 \begin{theorem}
 {There exists a function $F=(f(0)=0,f(1), \dots)\in \bQ T$ with the following
 property: for each $n$ and for each $F^n$-admissible vector
 $\{h(i)\}, i=1, \dots, n$ (where $F^n=(f(0), \dots, f(n))\in \bQ T_n$
 is the initial fragment of $F$ of length $n$), there exists $N$ such
 that $f(N+i)=h(i), i=1,\dots, n$.}
 \end{theorem}

 Indeed, such a function $F$ satisfies condition (c) above, i.e.,
 $F$ is a universal Toeplitz function; and it defines a rational metric
 on $\Bbb Z$ with required property.

 In order to prove Theorem 8, we will construct such a function $F$
 by induction: we will define it by a recursive procedure, each
 step being founded on the following theorem:

 \begin{theorem}
 {For any given finite function $F_n=(f(i))_{i=0}^n \in \bQ T_n$ (or ${\Bbb Z}T_n$)
 and any vector $H_n=(h(i))_{i=1}^n$ that is admissible
 with respect to $F_n$,
 there exist a positive integer $m$ and a vector $G =(g(1), \dots, g(m))$ such that
 the function  ${\bar F}=(f(0),f(1), \dots,
 f(n), g(1), \dots, g(m), h(1), \dots, h(n))$
 (the prolongation of $F_n$) belongs to $\bQ T_{2n+m}$
 (or $\bZ T_{2n+m}$).}
 \end{theorem}

 Using this fact, the required universal function $F\in \bQ T$
 (or $F \in \bZ T$) can be easily obtained as follows. Enumerate all
 finite integer vectors $V=\{v_1,v_2, \dots \}$; choose an initial function
 $F_0 \in \bQ T_{n_0}$ and then find the first vector from the sequence $V$ that
 is admissible for $F_0$ and apply the theorem giving the prolongation
 of $F_0$; if we already have some finite  $F_n$, we
 choose the next admissible vector from $V$ and apply the prolongation procedure.
 As a result, we obtain an infinite universal function $F$.

 It suffices to prove Theorem 9 above for an integer
 function $F_n \in \bZ T_n$, because before applying
 the prolongation procedure we can
 multiply the given rational vector $F \in \bQ T_n$ by a suitable integer and
 obtain $F_n \in \bZ T_n$, and then, after the prolongation procedure,
 divide the prolongation by the same integer.
 So we will give a proof of Theorem 8 for the case $F_n \in {\bZ}T_n$.

 The construction of prolongation is based on several lemmas,
 the first one being the simple {\it amalgamation lemma}, which is an analog of
 the lemma from \cite{vershik1} and \cite{hco}.

\begin{lemma}
 For each $F_n \in \bQ T_n$, let
 \[M(F_n)=\max_{1\leq k \leq n}|f(k)-f(n-k+1)|,\qquad
 m(F_n)=\min_{1\leq k \leq n}(f(k)+f(n-k+1).\] Then $M(F_n)\leq
 m(F_n)$.
\end{lemma}
\begin{pf}
 {We have $f(k)+f(n-k)+f(n-j)=d(0,k)+d(0,n-k)+d(n,j)=d(0,k)+d(n,k)+
 d(n,j) \geq d(0,k)+d(k,j)\geq d(0,j)=f(j)$, so $f(k)+f(n-k)\geq
 f(j)-f(n-j)$ for all $j,k=1,\dots, n$. Consequently, $$m(F_n) \equiv
 \min_k [f(k)+f(n-k)] \geq \max_j |f(j)-f(n-j)|\equiv M(F_n).$$}
\qed
\end{pf}

The next lemma shows how to extend a given Toeplitz function and
a given admissible vector with one coordinate each, namely, how to
join a new coordinate to the function and a new coordinate to the
beginning of the admissible vector. It will be the base of induction.

\begin{lemma}
{Suppose we have a function
$F\equiv F_n=(f(i))_{i=1}^n \in \bZ T_n$\footnote{We omit the coordinate $f(0)$,
which is identically equal to $0$.}
 and an $F_n$-admissible integer vector $H\equiv H_n=(h(i))_{i=1}^n$. Set
$$d=\max_{i=1,\dots, n}\{f_i,h_i\}.$$
Then there exist two integers $g(1),g(N)$ such that

{\rm(i)} the function $F_{n+1}=(f(1), \dots, f(n),g(1))$ lies in $\bZ T_{n+1}$;

{\rm(ii)} the vector $H_{n+1}=(g(N),h(1), \dots, h(n))$ is $F_{n+1}$-admissible;

{\rm (iii)} $2\leq g(1)\leq d-1$, $2\leq g(N) \leq d-1 $.}
\end{lemma}

\begin{pf}
{The numbers $g(1)$ and $g(N)$ are solutions of the following system of
linear inequalities:
\begin{eqnarray} &&\max_{i=1,\dots, n}|f(i)-f(n-i)|\leq g(1) \leq
\min_{i=1,\dots, n}(f(i)+f(n-i))\\
 &&\max_{i=1,\dots, n}|f(i)-h(i)|\leq g(N) \leq \min_{i=1,\dots, n}(f(i)+h(i)),
\\
&&|g_1-g_N|\leq h_n \leq (g(1)+g(N)),
\\
&& 2\leq g(1) \leq d-1, \quad 2 \leq g(N) \leq d-1.
\end{eqnarray}
 Inequality (1) expresses the fact that the extended function
 $F'\equiv (f(1),\dots,$ $f(n),g(1))$
 belongs to $\bZ T_n$; inequality (2) means that the
 extended vector $(g(N),h(1),\dots, h(n))$ is $F'$-admissible; inequality (3),
 together with inequality (1) means the $F'$-admissibility of the coordinate $h(n)$.
 The compatibility of inequalities (1)--(3) is a corollary of the amalgamation
 lemma
 and the $F$-admissibility of $h$, which is easy to check;
 condition (4) follows from the definition
 of $d$ and easy calculations. So the required properties
 (i), (ii), (iii) for this choice of $g(1)$ and $g(N)$ are fulfilled.}\qed
\end{pf}

The construction of the vector $(g_1, \dots, g_m)$ is as follows.
First of all, we put $m=nd$ and define a prolongation of $F$ as a
vector divided into $d$ blocks, each of length $n$:
$$(F,G_1, \dots, G_d, H) \mbox{ where } G_i = (g((i-1)n+1),
\dots, g(in));\quad i=1, \dots,d.$$ It is convenient to denote $F=G_0$ and
$H=G_{d+1}$; the coordinates $f(i)$, $h(i)$ belong to the integer
interval $[1,d]$.  We will define vectors $G_i$, $i=1, \dots, d$
successively in the following order (``from both sides''): $G_1,G_d$,
then $G_2,G_{d-1}$, and so on up to
$G_{\frac{d-3}{2}},G_{\frac{d+1}{2}}$, and, finally,
$G_\frac{d-1}{2}$. Now let us define $G_1$ and $G_2$.

We apply Lemma 11 and obtain $g(1)$ and $g(N)\equiv g(nd)$. The condition of the lemma is fulfilled
for the extended vector, and we can repeat the same procedure with the extended vectors
$F'=(f(1),\dots, f(n),g(1))$ and $H'=(g(nd),h(1), \dots, h(n))$ and join  numbers $g(2)$ and
$g(nd-1)$, obtaining vectors $(f(1),\dots,f(n),g(1),g(2))$ and $=(g(nd-1),g(nd),h(1), \dots,
h(n))$; then we join $g(3)$ and $g(nd-2)$, etc., up to $g(n)$ and $g((n-1)d+1)$. All these integers
belong to the interval $[2,d-1]$. As a result, we obtain the first part of the construction,
namely, the vectors $G_1, G_d$. We continue in the same way and apply Lemma 10 to the vectors
$G_1,G_2$ (instead of $F,H$) with only one change: the integer coordinates $g(n+1),g(n+2),\dots,
g(2n)$ and $g((n-1)d),\dots, g((n-2)d+1)$ belong to the interval $[3,d-2]$ (instead of $[2,d-1]$ in
the case of $G_1, G_d$). Thus the interval shrinks. By induction, we obtain all blocks up to the
last block $G_\frac{d-1}{2}$, whose all coordinates are equal:
$g(n(\frac{d-1}{2}+i))=\frac{d-1}{2}$ (because the interval is reduced to the single point
$\frac{d-1}{2}$). So we join the beginning and the end and obtain the required vector
$G_0=F,G_1,\dots, G_d,H$.\qed

This gives the proof of Theorems 8, 9 and the existence of a transitive isometry on ${\Bbb
Z}\bU,{\Bbb Q}\bU$ and an isometry on $\bU$ with dense orbit.

\bigskip

The proof of the theorem gives also further information:

\begin{cor}
The group $\Aut(\bQ\bU)$ contains $2^{\aleph_0}$ conjugacy classes of isometries
which permute the points in a single cycle. Moreover, representatives of
these classes remain non-conjugate in $\Aut(\bU)$.
\end{cor}

\begin{pf} It is clear that, if cyclic isometries $g$ and $h$ are
conjugate, then the functions $f_g$ and $f_h$ describing them as in the
above proof are equal. For, if $h=k^{-1}gk$, then
\[f_h(n)=d(x,h^n(x))=d(x,k^{-1}g^nk(x))=d(k(x),g^nk(x))=f_g(n).\]
But the set of functions describing cyclic isometries of $\bQ\bU$ is residual,
hence of cardinality~$2^{\aleph_0}$.\qed
\end{pf}

The cyclic isometries constructed in this section have the property that
$d(x,g(x))$ is constant for $x\in \bQ\bU$, and hence this holds for all
$x\in \bU$. In particular, these isometries are bounded.

\section{An abelian group of exponent~$2$}

To extend this argument to produce other groups acting regularly on~$\bQ\bU$,
it is necessary to change the definition of a Toeplitz function so that
the metric is defined by translation in the given group. We give here one
simple example.

\begin{prop}
The countable abelian group of exponent~$2$ acts regularly as an isometry
group of $\bQ\bU$.
\end{prop}

\begin{pf}
This group $G$ has a chain of subgroups $H_0\le H_1\le H_2\le\cdots$ whose
union is $G$, with $|H_i|=2^i$. We show that, given any $H_i$-invariant
rational metric on $H_i$ and any $h\in H_{i+1}\setminus H_i$, we can
prescribe the distances from $h$ to $H_i$ arbitrarily (subject to the
consistency condition) and extend the result to an $H_{i+1}$-invariant
metric on $H_{i+1}$. The extension of the metric is done by translation
in $H_{i+1}$: note that $H_{i+1}\setminus H_i$ is isometric to $H_i$,
since $d(h+h',h+h'')=d(h',h'')$ for $h',h''\in H_i$. Now the resulting
function is a metric. All that has to be verified is the triangle
inequality. Now triangles with all vertices in $H_i$, or all vertices in
$H_{i+1}\setminus H_i$, clearly satisfy the triangle inequality. Any
other triangle can be translated to a triangle containing $h$ and two
points of $H_i$, for which the triangle inequality is equivalent to the
consistency condition for extending the metric to $H_i\cup\{h\}$.\qed
\end{pf}

Note that almost all $G$-invariant metrics (in the sense of
Baire category) are isometric to $\bQ\bU$.

We can construct the analogous transitive actions for other
abelian groups on $\bU$.

\begin{prop}
There are transitive abelian groups of isometries of $\bU$ of the
groups of exponent~$2$.
\end{prop}

\begin{pf} Let $G$ be one of the abelian groups previously constructed,
and $\overline{G}$ its closure in $\Aut(\bU)$. Since the orbits of
$G$ are dense, it is clear that $\overline{G}$ is transitive.
Moreover, as the closure of an abelian group, it is itself
abelian. For, if $h,k\in\overline{G}$, say $h_i\to h$ and $k_i\to
k$; then $h_ik_i=k_ih_i\to hk=kh$. Similarly, if $G$ has
exponent~$2$, then so does $\overline{G}$.\qed
\end{pf}

\section{Other regular group actions}

There is a one-way relationship between transitive group actions on $Q$ (or, more generally, group
actions on $\bU$ with a dense orbit) and transitive actions on the Universal (random) graph $R$, as
given in the following result.

\begin{prop}
Let $G$ be a group acting on Urysohn space $\bU$ with a countable dense orbit $X$. Then there
exists a natural structure of the universal graph $R$ on~$X$ and group $G$ preserves this
structure.
\end{prop}

\begin{pf}
Partition the positive real numbers into two subsets $E$ and $N$ such that,
for any $R,\epsilon>0$, there are consecutive intervals of length at most
$\epsilon$ to the right of $R$ with one contained in $E$ and the other in $N$.
(For example, take a divergent series $(a_n)$ whose terms tend to zero, and
put half-open intervals of length $a_n$ alternately in $E$ and $N$.)

We define a graph on $X$ by letting $\{x,y\}$ be an edge if $d(x,y)\in E$,
and a non-edge if $d(x,y)\in N$. Clearly this graph is $G$-invariant; we must
show that it is isomorphic to the random graph $R$.

Let $U$ and $V$ be finite disjoint sets of points of $X$, and let the
diameter of $U\cup V$ be $d$ and the minimum distance between two of its
points be $m$. Choose $R>d/2$ and $\epsilon<m/2$, and find consecutive
intervals $I_E$ and $I_N$ as above. Let $U\cap V=\{w_1,\ldots,w_n\}$.
Choosing arbitrary values $g(i)\in I_E\cap I_n$, the consistency condition
\[|g(i)-g(j)|\le d(w_i,w_j)\le g(i)+g(j)\]
is always satisfied. So choose the values such that $g(i)$ is in
the interior of $I_E$ if $w_i\in U$, and in the interior of $I_n$
if $w_i\in V$. Let $z$ be a point of $U$ with $d(z,w_i)=g(i)$.
Since $X$ is dense, we can find $x\in X$ such that $d(x,z)$ is
arbitrarily small; in particular, so that $d(x,w_i)$ is in $I_E$
(resp. $I_N$) if and only if $d(z,w_i)$ is. Thus $x$ is joined to
all vertices in $U$ and to none in $V$. This condition
characterizes $R$ as a countable graph.\qed
\end{pf}

The converse is not true. A special case of the result of Cameron and
Johnson~\cite{cj} shows that a sufficient condition for a group $G$ to act
regularly on $R$ is that any element has only finitely many square roots. In
a group with odd exponent, each element has a unique square root. So any such
group acts regularly on $R$. But we have the following:

\begin{prop}
The countable abelian group of exponent~$3$ cannot act on $\bU$ with a dense
orbit, and in particular cannot act transitively on $\bQ\bU$.
\end{prop}

\begin{pf} Suppose that we have such an action of this group $A$. Since the
stabiliser of a point in the dense orbit is trivial, we can identify the
points of the orbit with elements of $A$ (which we write additively).

Choose $x\ne 0$ and let $d(0,x)=\alpha$. Then $\{0,x,-x\}$ is an equilateral
triangle with side~$\alpha$. Since $\bU$ is universal and $A$ is dense, there
is an element $y$ such that $d(x,y),d(-x,y)\approx\frac{1}{2}\alpha$ and
$d(0,y)\approx\frac{3}{2}\alpha$. (The approximation is to within a given
$\epsilon$ chosen smaller than $\frac{1}{6}$. Then the three points $0,y,x-y$
form a triangle with sides approximately $\frac{3}{2}\alpha$,
$\frac{1}{2}\alpha$, $\frac{1}{2}\alpha$, contradicting the triangle
inequality.\qed
\end{pf}

\section{Unbounded isometries of $\bU$}

The subgroup $\BA(\bU)$ is not the whole isometry group, because unbounded isometries exist. The
simplest way to see this is to mention that  Euclidean space $R^n$ can be imbedded to $\bU$ in such
a way that the group of motions $Iso(R^n$ has monomorphic imbedding to $Iso(\bU)$. But we will give
a direct construction of such isometry. (We are grateful to J.Ne\v{s}et\v{r}il for the following
argument.)

\begin{prop}
There exist unbounded isometries of~$\bQ\bU$ (and hence of $\bU$).
\end{prop}

The proof depends on a lemma.

\begin{lemma}
Let $A$ be a finite subset of $\bQ\bU$ and let $g$ be a function on $A$
satisfying  the consistency conditions ($*$). Then the diameter of the set
\[S_A=\{z\in \bQ\bU: d(z,a)=g(a)\hbox{ for all }a\in A\}\]
is twice the minimum value of $g$.
\end{lemma}

\begin{pf} Let $z_1$ and $z_2$ be two points from $S_A$. Consider
the problem of adding $z_2$ to the set $A\cup\{z_1\}$. The consistency
conditions for $z_2$ are precisely those for $z_1$ together with the
conditions
\[|d(z_2,z_1)-d(z_2,a)|\le d(z_1,a)\le d(z_2,z_1)+d(z_2,a)\]
for all $a\in A$. Since $d(z_1,a)=d(z_2,a)=g(a)$, the only non-trivial
restriction is $d(z_1,z_2)\le 2d(z_1,a)=2g(a)$, which must hold for all
$a\in A$.\qed
\end{pf}

\paragraph{Proof of Proposition 17} We construct an isometry $f$ of $\bQ\bU$ by
the standard back-and-forth method, starting with any enumeration of $\bQ\bU$.
At odd-numbered stages we choose the first point not in the range of $f$
and select a suitable pre-image. At stages divisible by~$4$ we choose the
first point not in the domain of $f$ and select a suitable image. This
guarantees that the isometry we construct is a bijection from $\bQ\bU$ to
itself.

At stage $4n+2$, let $\bU$ be the domain of $f$. Choose an unused point $z$
whose least distance from~$\bU$ is $n$. Now the diameter
of the set of possible images of $z$ is $2n$; so we can choose a
possible image $f(z)$ whose distance from~$z$ is at least $n$. Then the
constructed isometry is not bounded.\qed

\bigbreak

We can improve this argument to construct an isometry $g$ such that all
powers of $g$ except the identity are unbounded. In fact, even more is true:

\begin{lemma}
There are two isometries $a,b$ of $\bQ\bU$ which generate a free group, all of
whose non-identity elements are unbounded isometries.
\end{lemma}

\begin{pf}
We begin by enumerating $\bQ\bU=(x_0,x_1,\ldots)$. We follow the argument we used
to show that unbounded isometries exist. We construct $a$ and $b$
simultaneously, using the even-numbered stages for a back-and-forth argument
to ensure that both are bijections, and the odd-numbered stages to ensure
that any word in $a$ and $b$ is unbounded. The first requirement is done as
we have seen before.

Enumerate the words $w(a,b)$ in $a$ and $b$ and their inverses. (It suffices
to deal with the cyclically reduced words, since all others are conjugates of
these.) We show first how to ensure that $w(a,b)\ne1$. At a given stage,
suppose we are considering a word $w(a,b)$. Choose a point $x_i$ such that
neither $a$ nor $b$, nor their inverses, has been defined on $x_j$ for $j\ge
i$. Suppose that $w$ ends with the letter $a$. Since there are infinitely
many choices for the image of $x_i$ under $a$, we may choose an image $x_j$
with $j>i$. Now define the action of the second-last letter of the word on
$x_j$ so that the image is $x_k$ with $k>j$. Continuing in this way, we end
up with a situation where $w(a,b)x_i=x_m$ with $m>i$. So $w(a,b)\ne1$.

To ensure that $w(a,b)$ is unbounded, we must do more. Enumerate the words so
that each occurs infinitely often in the list. Now, the $k$th time we revisit
the word $w$, we can ensure (as in our construction of an unbounded isometry)
that $d(x_i,w(a,b)x_i)\ge k$. Thus $w(a,b)$ is unbounded.\qed
\end{pf}

\section{A dense free subgroup of $\Aut(\bU)$}

We can now use a trick due to Tits~\cite{bldgs} to show that there is a dense
subgroup of $\Aut(\bU)$ which is a free group of countable rank.

\begin{theorem}
There is a subgroup $F$ of $\Aut(\bQ\bU)$ which acts faithfully and homogeneously
on $\bQ\bU$ and is isomorphic to the free group of countable rank.
\end{theorem}

\begin{pf}
Since the free group $F_2$ contains a subgroup isomorphic to $F_\omega$,
choose a group $H$ with free generators $h_i$ for $i\in\mathbb{N}$, such that
$H\cap\BA(\bQ\bU)=1$. Enumerate the pairs of isometric $n$-tuples of elements of
$\bQ\bU$, for all $n$, as $(\alpha_0,\beta_0)$, $(\alpha_1,\beta_1)$, \dots~.
Now, for each $i$, Lemma~\ref{bddense} shows that we can choose
$n_i\in\BA(\bQ\bU)$ such that $n_ih_i(\alpha_i)=\beta_i$. Let $F$ be the group
generated by the elements $n_1h_1, n_2h_2, \ldots$~. Clearly $F$ acts
homogeneously on $\bQ\bU$. We claim that $F$ is free with the given generators.
Suppose that $w(n_ih_i)=1$ for some word $w$. Since $\BA(\bQ\bU)$ is a normal
subgroup, we have $nw(h_i)=1$ for some $n\in\BA(\bQ\bU)$. Since $n$ is bounded
and $w(h_i)$ unbounded, this is impossible. In fact this argument shows that
all the non-identity elements of $F$ are unbounded isometries.\qed
\end{pf}

\end{document}